\documentclass[12pt]{article}
\usepackage{graphicx}
\usepackage{amsmath}
\usepackage{amsfonts}
\usepackage{amssymb}
\usepackage{amscd}
\usepackage{amsthm}

\theoremstyle{plain}
\newtheorem{theorem}{Theorem}
\newtheorem{corollary}[theorem]{Corollary}

\newtheorem{proposition}[theorem]{Proposition}

\theoremstyle{definition}

\newtheorem{remark}[theorem]{Remark}

\def\o{\omega}
\def\a{\alpha}
\def\aa{\eta}
\def\b{\beta}
\def\bb{\gamma}

\def\tr{\mathbf{t}}
\def\t{t}

\def\ps{\varphi}
\def\f{\psi}

\def\A{\mathcal{A}}
\def\B{\mathcal{B}}

\def\S{\mathbf{S}}
\def\SS{\mathbf{S}}

\def\E{\mathcal{E}}
\def\K{\mathcal{K}}
\def\R{\mathbb{R}}

\def\Z{\mathbb{Z}}

\date{}

\begin{document}

\title{All strongly-cyclic branched coverings of $(1,1)$-knots are Dunwoody manifolds}

\author{Alessia Cattabriga \and Michele Mulazzani}

\maketitle


\begin{abstract}
We show that every strongly-cyclic branched covering of a
$(1,1)$-knot is a Dunwoody manifold. This result, together with
the converse statement previously obtained by Grasselli and
Mulazzani, proves that the class of Dunwoody manifolds coincides
with the class of strongly-cyclic branched coverings of
$(1,1)$-knots. As a consequence, we obtain a parametrization of
$(1,1)$-knots by 4-tuples of integers. Moreover, using a
representation of $(1,1)$-knots by the mapping class group of the
twice punctured torus, we provide an algorithm which gives
the parametrization  of all  torus knots in $\SS^3$.\\
\\ {{\it Mathematics Subject
Classification 2000:} Primary 57M12, 57N10; Secondary  57M25.\\
{\it Keywords:} $(1,1)$-knots, cyclic branched coverings, Heegaard
diagrams, Dunwoody manifolds, torus knots, cyclically presented
groups.}

\end{abstract}


\section{Introduction} \label{intro}

In order to investigate the relations between cyclic branched
coverings of knots in $\SS^3$ and manifolds admitting cyclically
presented fundamental groups, M. J. Dunwoody introduced in
\cite{D} a class of 3-manifolds depending on six integer
parameters. As proved in \cite{GM}, all these manifolds turn out
to be strongly-cyclic coverings of lens spaces (possibly $\S^3$),
branched over $(1,1)$-knots. Moreover, it has been shown in
\cite{Mul} that every $n$-fold strongly-cyclic branched covering
of a $(1,1)$-knot admits a genus $n$ Heegaard diagram encoding a
cyclic presentation for the fundamental group. This result has
been improved in \cite{CM}, obtaining a constructive algorithm
which, starting from a representation of $(1,1)$-knots through the
elements of the mapping class group of the twice punctured torus,
explicitly gives the cyclic presentations.

In this paper we prove that all strongly-cyclic branched coverings
of $(1,1)$-knots are actually Dunwoody manifolds. As a
consequence, the class of Dunwoody manifolds coincides with the
class of strongly-cyclic branched coverings of $(1,1)$-knots.

We also obtain, as a further consequence, a parametrization of all
$(1,1)$-knots (with the exception of the ``core'' knot
$\{P\}\times\S^1\subset\S^2\times\S^1$, which admits no
strongly-cyclic branched coverings) by means of four of the six
Dunwoody parameters. Moreover, we give an algorithm that allows us
to find the parametrization of all torus knots in $\SS^3$.

We refer to \cite{Ka,BZ} for details on knot theory and cyclic
branched coverings of knots, and to \cite{Jo} for details on
cyclic presentations of groups.

\section{Strongly-cyclic branched coverings of (1,1)-knots and Dunwoody manifolds}

An $n$-fold cyclic covering of a 3-manifold $N^3$ branched over a
knot $K\subset N^3$ is called {\it strongly-cyclic\/} if the
branching index of $K$ is $n$ (i.e., the fiber of each point of
$K$ contains a single point). So the homology class of a meridian
loop $m$ around $K$ is mapped by the associated monodromy
$\o:H_1(N^3-K)\to\Z_n$ to a generator of $\Z_n$ (up to equivalence
we can always suppose $\o[m]=1$).

Observe that a cyclic branched covering of a knot $K$ in $\S^3$ is
always strongly-cyclic and uniquely determined, up to equivalence,
since \hbox{$H_1(\S^3-K)\cong\Z$.} Obviously, this property is no
longer true for a knot in a more general 3-manifold. Also, if $p$
is a prime number, any $p$-fold cyclic branched covering of a knot
$K$ is automatically strongly-cyclic.

In this paper we deal with strongly-cyclic branched coverings of
$(1,1)$-knots, which are knots in lens spaces (possibly in
$\S^3$).

A knot $K$ in a 3-manifold $N^3$ is called a $(1,1)$-{\it knot\/}
if there exists a Heegaard splitting of genus one
$$(N^3,K)=(H,A)\cup_{\ps}(H',A'),$$ where $H$ and $H'$ are solid
tori, $A\subset H$ and $A'\subset H'$ are properly embedded
trivial arcs, and $\ps:(\partial H',\partial A')\to(\partial
H,\partial A)$ is an attaching homeomorphism (see Figure \ref{Fig.
1}). Obviously, $N^3$ turns out to be a lens space $L(p,q)$
(including $\S^3=L(1,0)$).


\begin{figure}[ht]
\begin{center}
\includegraphics*[totalheight=5.5cm]{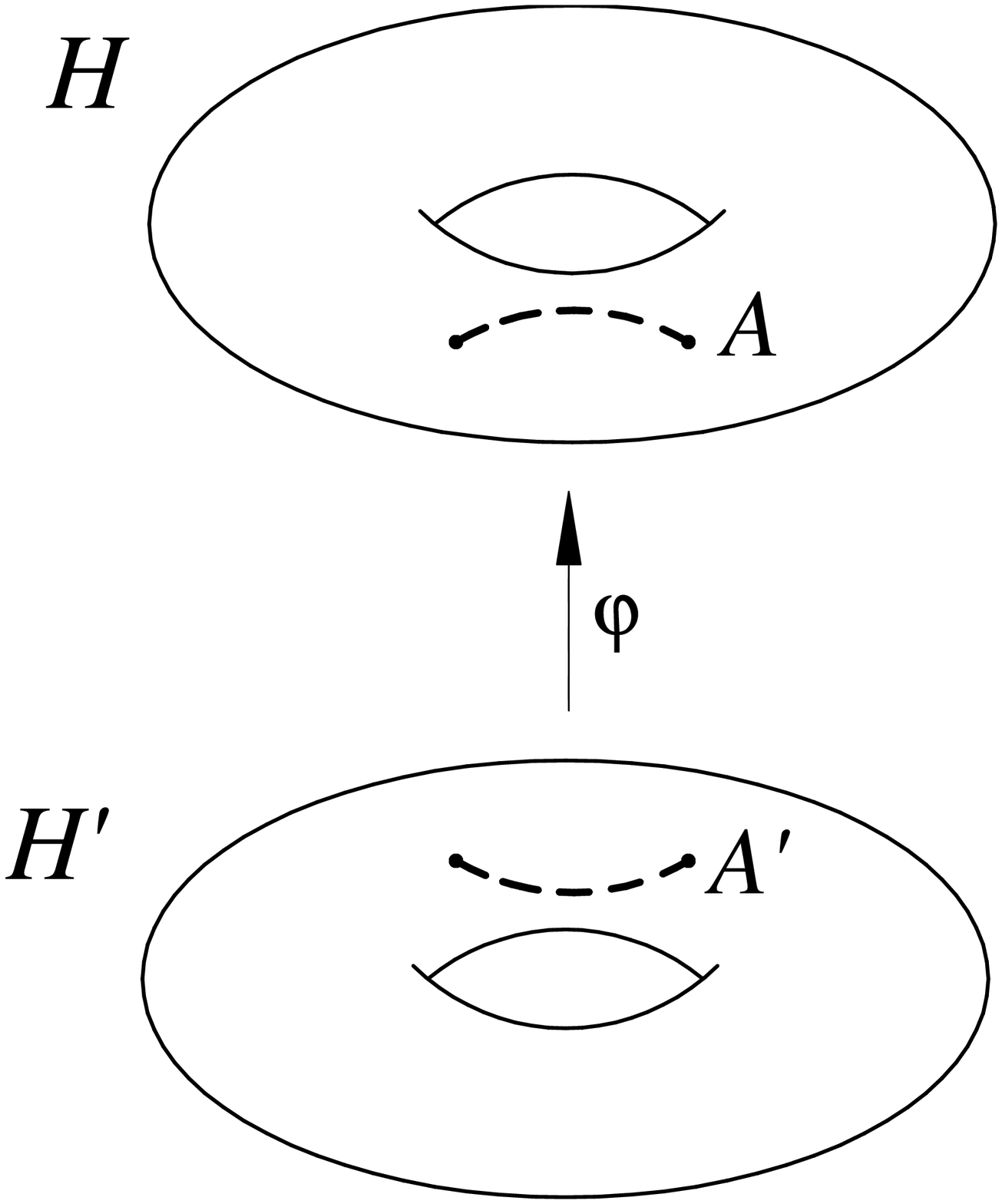}
\end{center}
\caption{A $(1,1)$-decomposition.} \label{Fig. 1}
\end{figure}

It is well known that the family of $(1,1)$-knots contains all
torus knots and all two-bridge knots in $\S^3$. Several
topological properties of $(1,1)$-knots have recently been
investigated (see references in \cite{CM2}).

\begin{proposition} \label{basic}
A $(1,1)$-knot $K\subset L(p,q)$ with $(1,1)$-decomposition
$(L(p,q),K)=(H,A)\cup_{\ps}(H',A')$ is completely determined, up
to equivalence, by $\ps(\beta')$, where $\beta'$ is the boundary
of a meridian disk $D'\subset H'$ which does not intersect $A'$.
Moreover, if $(L(p,q),\bar K)=(H,A)\cup_{\bar{\ps}}(H',A')$ is a
decomposition of a $(1,1)$-knot $\bar K$ such that
$\bar{\ps}(\beta')$ is isotopic to $\ps(\beta')$ in $\partial
H-\partial A$, then $\bar K$ is equivalent to $K$.
\end{proposition}
\begin{proof}
The first statement follows from the fact that two properly
embedded trivial arcs in a ball $B$, with the same endpoints,  are
isotopic rel $\partial B$. The second statement is
straightforward.
\end{proof}

An algebraic representation of $(1,1)$-knots has been developed in
\cite{CM} and \cite{CM2}, where it is shown that there is a
natural surjective map $$\f\in PMCG_2(\partial H)\mapsto K_{\f}\in
\K_{1,1}$$ from the pure mapping class group of the twice
punctured torus $PMCG_2(\partial H)$ to the class $\K_{1,1}$ of
all $(1,1)$-knots. Using this representation, the necessary and
sufficient conditions for the existence and uniqueness of an
$n$-fold strongly-cyclic branched covering of a $(1,1)$-knot have
been obtained (see \cite{CM}).

\medskip

The family of Dunwoody manifolds has been introduced in \cite{D}
by a class of trivalent regular planar graphs (called {\it
Dunwoody diagrams\/}), depending on six integers $a,b,c,n,r,s$,
such that $n>0$, $a,b,c\ge 0$. For certain values of the
parameters, called {\it admissible\/}, the Dunwoody diagrams
$D(a,b,c,n,r,s)$ turn out to be Heegaard diagrams, hence defining
a wide class of closed, orientable 3-manifolds $M(a,b,c,n,r,s)$
with cyclically presented fundamental groups, called {\it Dunwoody
manifolds}.

\begin{figure}
 \begin{center}
 \includegraphics*[totalheight=6.5cm]{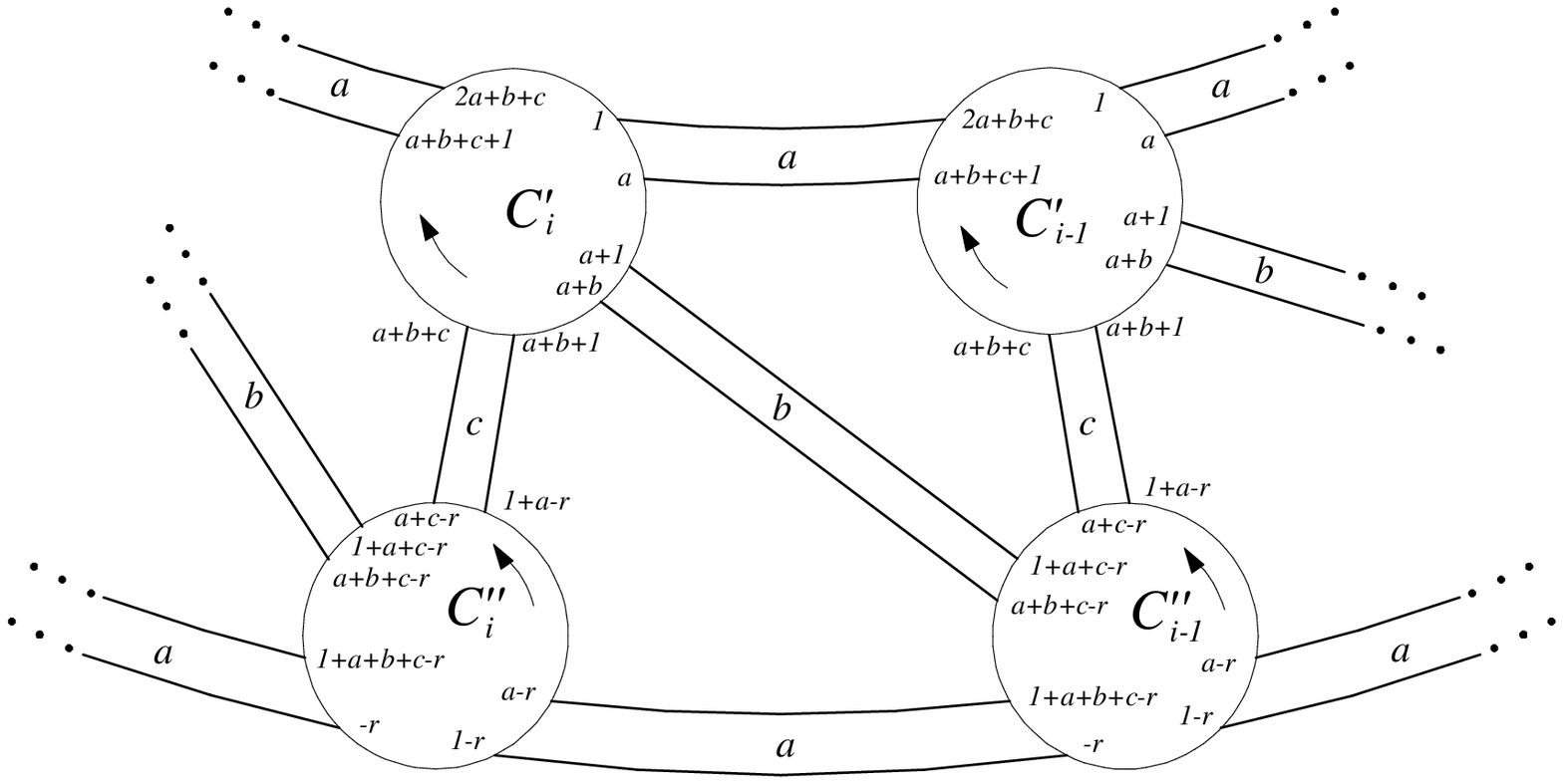}
 \end{center}
 \caption{The diagram $D(a,b,c,n,r,s)$, for $a+b+c>0$.}
 \label{Fig. 13}
\end{figure}

\begin{figure}
 \begin{center}
 \includegraphics*[totalheight=6.5cm]{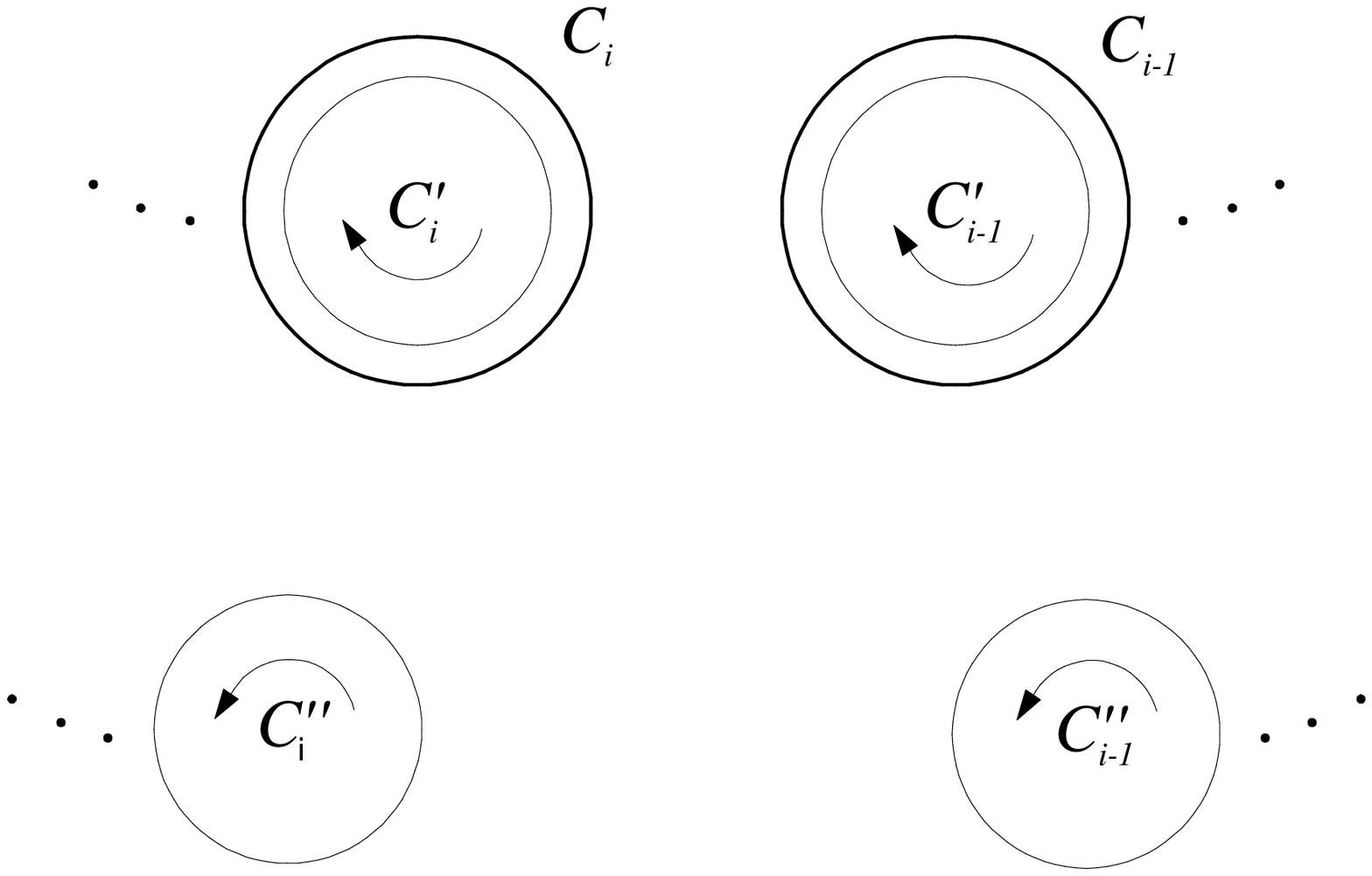}
 \end{center}
 \caption{The diagram $D(0,0,0,n,r,s)$.}
 \label{Fig. 28}
\end{figure}

More precisely, an admissible Dunwoody diagram $D(a,b,c,n,r,s)$ is
an open Heegaard diagram of genus $n$, with cyclic symmetry of
order $n$. It contains $n$ internal circles $C'_1,\ldots,C'_n$,
and $n$ external circles $C''_1,\ldots,C''_n$, each having
$d=2a+b+c$ vertices. These circles represent the first system of
curves of the Heegaard splitting. If $d>0$, as shown in Figure
\ref{Fig. 13}, the circle $C'_i$ (resp. $C''_i$) is connected to
the circle $C'_{i+1}$ (resp. $C''_{i+1}$) by $a$ parallel arcs, to
the circle $C''_{i}$ by $c$ parallel arcs and to the circle
$C''_{i-1}$ by $b$ parallel arcs, for every $i=1,\ldots,n$
(subscripts mod $n$). If $d=0$ (i.e., $a=b=c=0$), there are no
arcs connecting the circles, and the diagram (called {\it
trivial\/}) contains other $n$ circles $C_1,\ldots,C_n$, as
depicted in Figure \ref{Fig. 28}.

We denote by $\E$ the set of arcs when $d>0$, or the set of curves
$C_1,\ldots,C_n$ when $d=0$. Obviously, $\E$ represents the second
system of curves of the Heegaard splitting. To reconstruct the
splitting, the circle $C'_i$ must be glued to the circle
$C''_{i+s}$, so that, when $d>0$, equally labelled vertices are
identified together. Observe that the parameters $r$ and $s$ can
be considered mod $d$ and $n$ respectively, and we can suppose
$r=0$ when $d=0$. Since the identification rule and the diagram
are invariant with respect to an obvious cyclic action of order
$n$, the Dunwoody manifold $M(a,b,c,r,n,s)$ admits a cyclic
symmetry of order $n$. Of course, $M(a,b,c,1,r,0)$ is homeomorphic
to a lens space or to $\S^3$, since it admits a genus one Heegaard
splitting. Moreover, the trivial case $M(0,0,0,n,0,s)$ is
homeomorphic to the connected sum of $n$ copies of
$\SS^2\times\SS^1$, for all $n$ and $s$.

A characterization of all Dunwoody manifolds as strongly-cyclic
branched coverings of $(1,1)$-knots is given by the following
result.

\begin{proposition} \label{strongly-cyclic} \textup{\cite{GM}}
The Dunwoody manifold $M(a,b,c,n,r,s)$ is the $n$-fold
strongly-cyclic covering of the lens space $M(a,b,c,1,r,0)$
(possibly $\S^3$), branched over a $(1,1)$-knot only depending on
the integers $a,b,c,r$.
\end{proposition}

An interesting example of a Dunwoody manifold is $M(1,1,1,3,2,1)$,
which is homeomorphic to $\S^1\times\S^1\times\S^1$. It is well
known that this manifold cannot be a cyclic branched covering of
any knot in $\S^3$, but turns out to be a 3-fold cyclic covering
of $\S^2\times\S^1\cong M(1,1,1,1,2,0)$, branched over a
$(1,1)$-knot, which will be referred to as $K(1,1,1,2)$.

In the next section we prove the converse of Proposition
\ref{strongly-cyclic}. As a consequence, the class of Dunwoody
manifolds coincides with the class of strongly-cyclic branched
coverings of $(1,1)$-knots.

\section{Main result}

Now we establish the main result of this paper.

\begin{theorem}
\label{ultimo} Every strongly-cyclic branched covering of a
$(1,1)$-knot is a Dunwoody manifold.
\end{theorem}
\begin{proof}
Let $K\subset L(p,q)$ be a $(1,1)$-knot and let
$(L(p,q),K)=(H,A)\cup_{\ps}(H',A')$ be a $(1,1)$-decomposition of
$K$. Let $\b$ (resp. $\b'$) be a meridian of $\partial H$ (resp.
$\partial H'$) that bounds a disc in $H$ (resp. $H'$) not
intersecting $A$ (resp. $A'$).  The system of curves
$(\b,\ps(\b'))$  on $T=\partial H$ defines a genus one Heegaard
diagram of $L(p,q)$, which does not intersect $\partial
A=\{N,S\}$. Let $H_{\ps}$ be the open Heegaard diagram on $\R^2$
obtained by cutting $T$ along $\b$, and considering $S$ as the
point at the infinity of $\S^2=\R^2\cup\{S\}$. The diagram
consists of two canonical circles $C'$ and $C''$, corresponding to
$\b$, and  a closed curve or a set of arcs with endpoints on the
canonical circles, which corresponds to $\ps(\b')$ and will be
denoted by $\E$. Suppose that one of the following holds:
\begin{itemize}
\item[(1)] $H_{\ps}$ is the diagram depicted in Figure
\ref{Fig27} a); \item[(2)] $H_{\ps}$ is the diagram depicted in
Figure \ref{Fig27} b); \item[(3)] there  exist integers $a,b,c$,
with $a,b,c\geq 0$ and $a+b+c>0$, such that $H_{\ps}$ is the
diagram depicted in Figure \ref{Fig26}.
\end{itemize}

\begin{figure}
\begin{center}
\includegraphics*[totalheight=6cm]{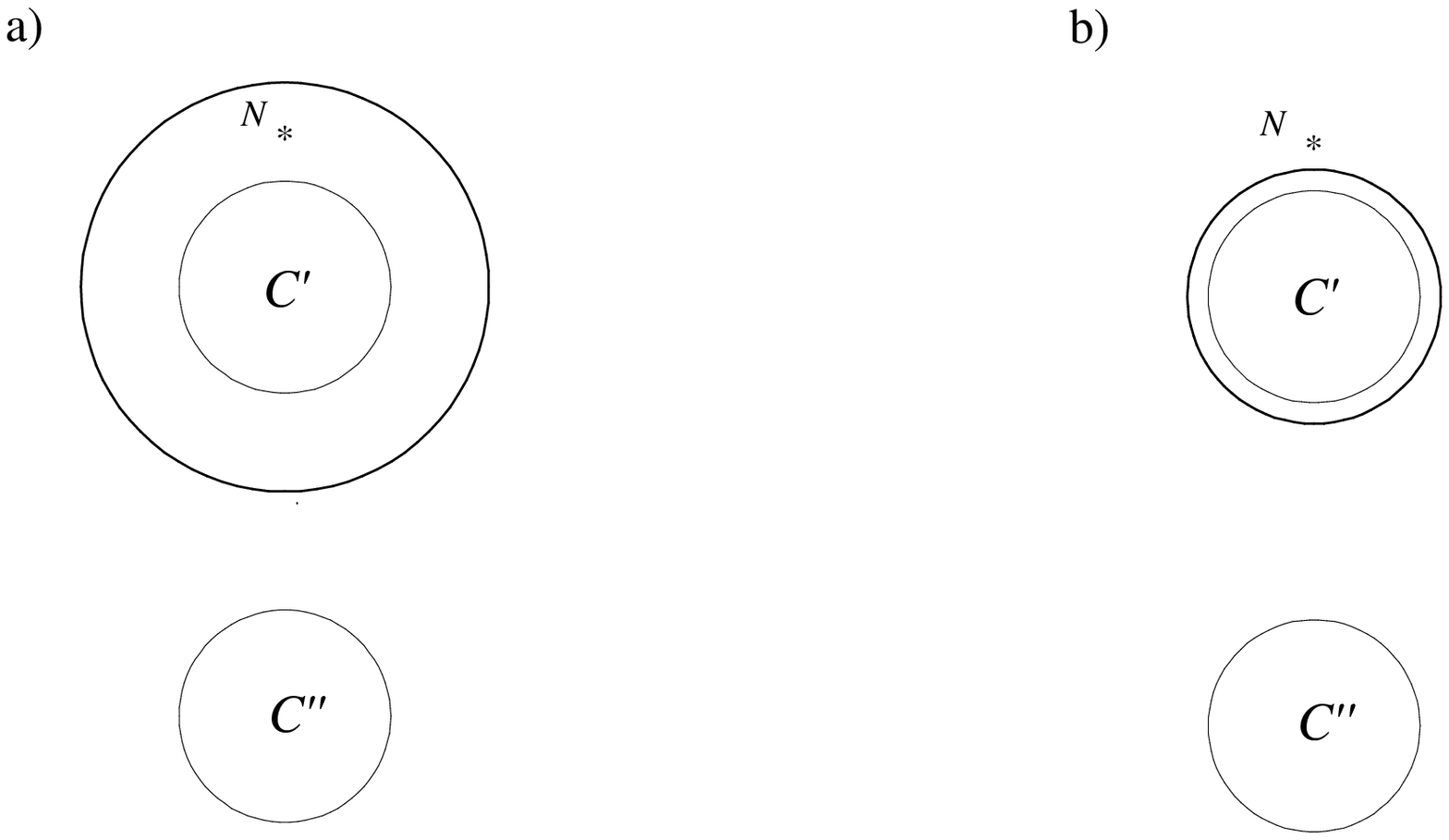}
\end{center}
\caption{} \label{Fig27}
\end{figure}

\begin{figure}
\begin{center}
\includegraphics*[totalheight=7cm]{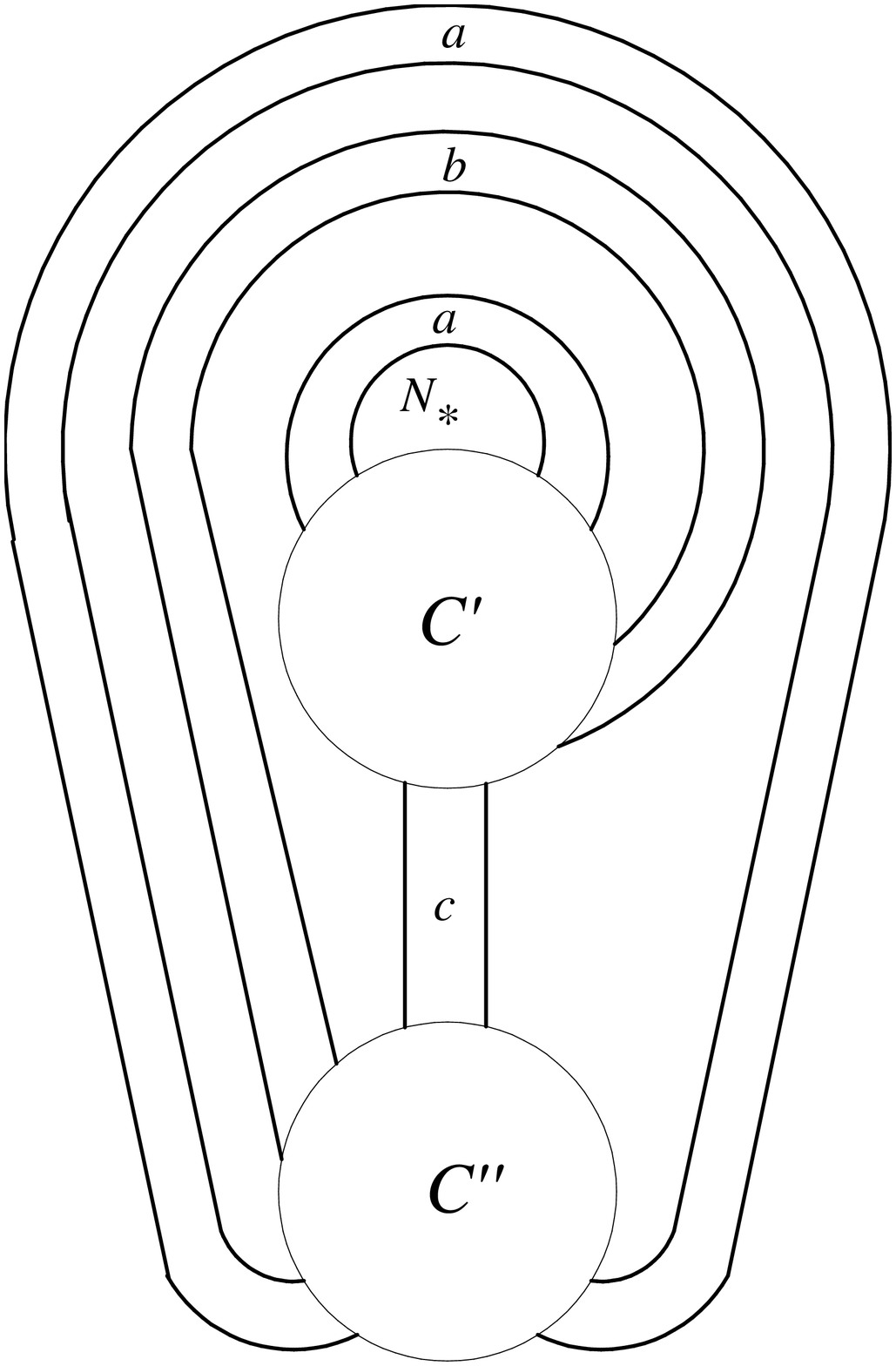}
\end{center}
\caption{} \label{Fig26}
\end{figure}

In the first case, $K$ is the core knot $\{P\}\times\SS^1\subset
\SS^2\times\SS^1$, where $P$ is a point of $\SS^2$. Therefore,
from \cite[Cor. 2]{CM}, we have $H_1(\SS^2\times\SS^1-K)=\langle
\a,\bb\,|\,\bb\rangle\cong\Z,$ where $\alpha$ and $\gamma$ are the
curves on $T$ depicted in Figure \ref{Fig12}. So, by \cite[Th.
4]{CM}, there exists no strongly-cyclic branched covering of $K$.

\begin{figure}[h]
\begin{center}
\includegraphics*[totalheight=4cm]{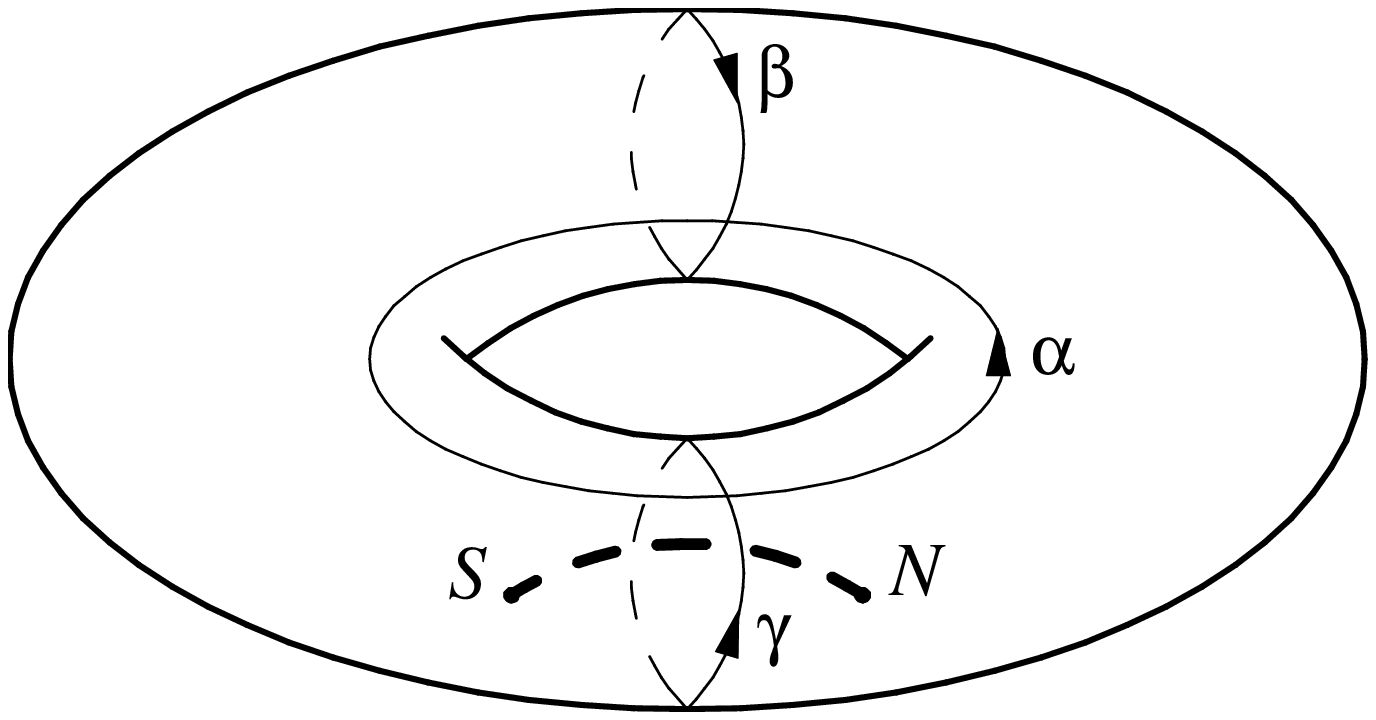}
\end{center}
\caption{} \label{Fig12}
\end{figure}

In the second case, $K$ is the trivial knot in $\SS^2\times\SS^1$.
Therefore, by \hbox{\cite[Cor. 2]{CM}}, we have
$H_1(\SS^2\times\SS^1-K)=\langle
\a,\bb\,|\,\emptyset\rangle\cong\Z\oplus\Z.$ So, by
\hbox{\cite[Th. 4]{CM},} there exist exactly $n$ $n$-fold
strongly-cyclic branched coverings of $K$, depending on the choice
of $\omega(\a)\in\Z_n$, where
$\omega:H_1(\SS^2\times\SS^1-K)\to\Z_n$ is the monodromy map of
the covering such that $\omega(\bb)=1$. If we denote by
$C_{n,s}(K)$ the $n$-fold strongly-cyclic branched covering of $K$
such that $\omega(\a)=s$, we have $C_{n,s}(K)=M(0,0,0,n,0,s)$.
Actually, as previously observed, $C_{n,s}(K)$ is homeomorphic to
the connected sum of $n$ copies of $\SS^2\times\SS^1$, for all
$n,s$.

Let us consider the third case. If $f:M\to L(p,q)$ is an $n$-fold
strongly-cyclic branched covering of $K$, then the
$(1,1)$-decomposition of $K$ lifts to a genus $n$ Heegaard
splitting for $M$ (see \cite{Mul}). Since $\omega(\bb)=1$, up to
equivalence, then the lifting of  $H_{\ps}$ is the Dunwoody
diagram $D(a,b,c,n,r,s)$, where $s=\omega(\a)$. In other words,
$M$ is the Dunwoody manifold $M(a,b,c,n,r,s)$.

By Proposition \ref{basic}, to prove the theorem it is enough to
show that $H_{\ps}$ is equivalent, up to Singer moves fixing $N$,
to one of the three diagrams discussed above.

Denote by $D'$ and $D''$  the disks of $\R^2$ bounded by $C'$ and
$C''$, respectively. Moreover, let $\A'$ (resp. $\A''$) be the set
of arcs of $\E$ with both the endpoints on $C'$ (resp. $C''$), and
denote by $\B$  the remaining arcs of $\E$. Of course,
$\vert\A'\vert=\vert\A''\vert$. An arc $e\in\A'$ (resp. $\A''$) is
called {\it trivial\/} if the closed curve $e\cup e'$, where $e'$
is one of the two arcs of $C'$ (resp. $C''$) with the same
endpoints of $e$, bounds a disc containing neither $N$ nor $D''$
(resp. $D'$). As illustrated in Figure \ref{Fig30}, each trivial
arc can be removed by a Singer move of type IB (see \cite{Si}).
So, up to equivalence, we can suppose that $H_{\ps}$ contains no
trivial arcs. Observe that this assumption implies that $e\cup e'$
bounds a disc in $\R^2$ containing the point $N$, for every
$e\in\A'\cup\A''$. In fact, if there exists a non trivial arc $e$
of $\A'$ (resp. of $\A''$) such that $e\cup e'$ bounds a disk $D$
in $\R^2$ not containing $N$, then $D$ contains $D''$ (resp. $D'$)
and therefore there exists a trivial arc in $\A''$ (resp. $\A'$).

\begin{figure}[h]
\begin{center}
\includegraphics*[totalheight=7cm]{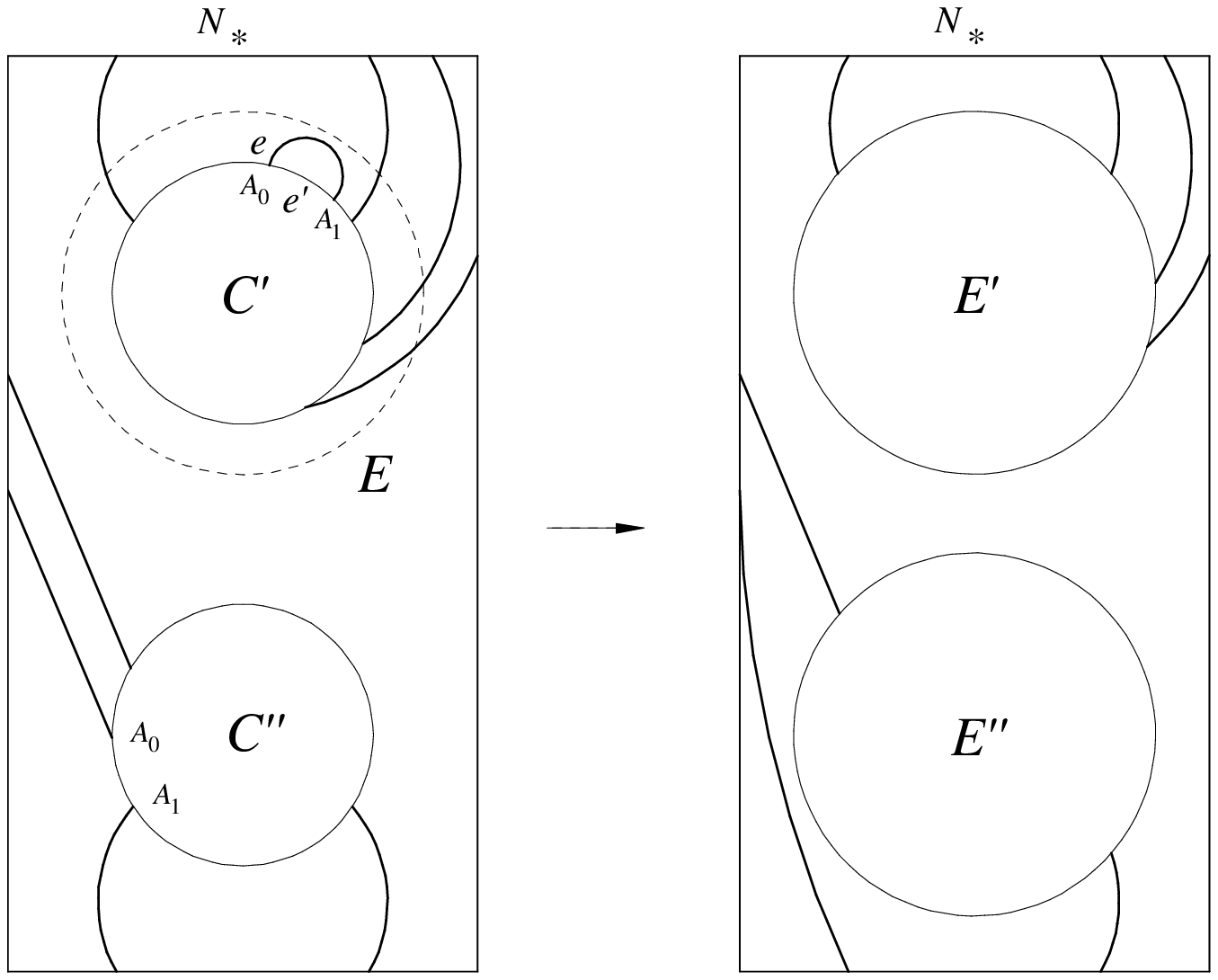}
\end{center}
\caption{Singer move of type  IB.} \label{Fig30}
\end{figure}

In order to simplify the proof, let us consider the planar graph
$\Gamma$ obtained from $H_{\ps}$ by collapsing the disks $D'$ and
$D''$ to their centers, that we still indicate by $C'$ and $C''$,
respectively. Of course, the arcs of $\A'$ and $\A''$ become loops
in $\Gamma$ bounding disks all containing $N$.

We say that two elements of $\E$ are {\it parallel\/} if they are
isotopic rel $\{C',C'',N\}$. It is easy to see that any two
elements of $\A'$ (resp. of $\A''$) are parallel. In fact, if the
disk bounded by a loop of $\A'$ (resp. $\A''$) contains $C''$
(resp. $C'$), then all the disks bounded by the loops of $\A'$
(resp. $\A''$) contain $C''$ (resp. $C'$). Otherwise, each loop of
$\A''$ (resp. $\A'$) bounds a disk not containing $N$. As regards
the elements of $\B$, we note that two different arcs $g,g'\in\B$
are parallel if and only if the closed curve $g\cup g'$ bounds a
disc $D_{g,g'}$ not containing $N$. It is not difficult to see
that there are at most two isotopy classes. For, if
$g,g',g''\in\B$ are different arcs such that $g$ is not parallel
to either $g'$ or $g''$, then $N\in D_{g,g'}$ and $N\in
D_{g,g''}$. Moreover, either $D_{g',g''}=(D_{g,g'}-D_{g,g''})\cup
g''$ or $D_{g',g''}=(D_{g,g''}-D_{g,g'})\cup g'$. In both cases
$N\notin D_{g',g''}$ and therefore $g'$ is parallel to $g''$.

\begin{figure}[h]
\begin{center}
\includegraphics*[totalheight=4.5cm]{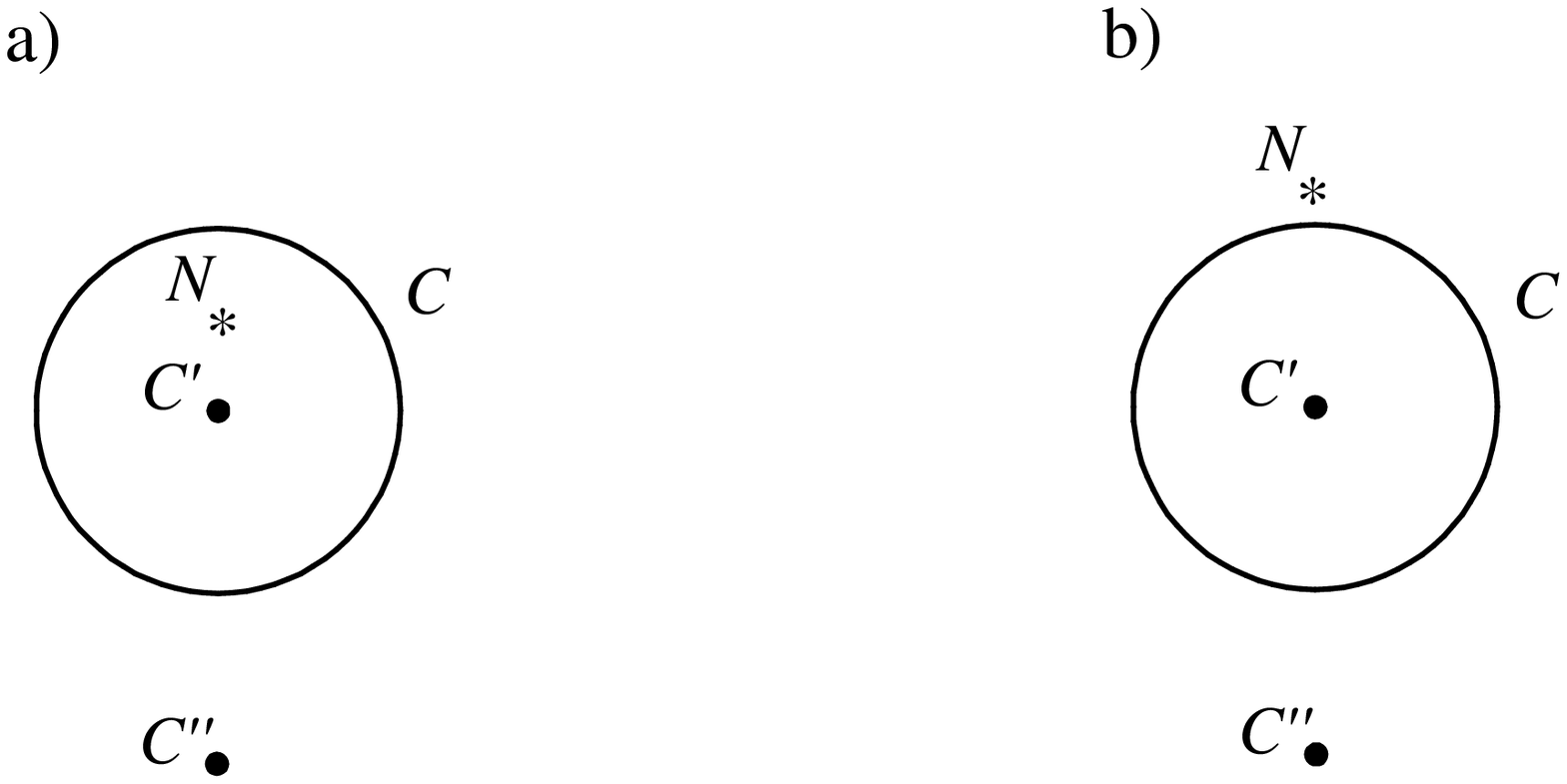}
\end{center}
\caption{} \label{Fig2}
\end{figure}

If $\A'=\A''=\B=\emptyset$, $\E$ consists of a  closed curve $C$.
So, up to isotopy in $\R^2-N$, we can suppose that $C$ is a
standard circle. There are two possibilities, depending on whether
the point $N$ is contained inside or outside $C$. But, in both
cases, since $C$ is a curve of a Heegaard diagram, $C'$ is inside
$C$ if and only if $C''$ is outside $C$. So, up to a possible
exchange between $C'$ and $C''$, the two possibilities are those
depicted in Figure \ref{Fig2}, which are the same as in Figure
\ref{Fig27}.

If $\A'\cup\A''\cup\B\ne\emptyset$, we can consider the graph
$\Gamma'$ obtained from $\Gamma$ by taking only one element for
each isotopy class of arcs. So $\Gamma'$ is a graph embedded in
$\R^2-N$ with two vertices, a loop in each vertex if
$\A'\ne\emptyset$, and one or two edges linking the vertices if
$\B\ne\emptyset$. If $\A'\ne\emptyset$, one of the two loops is
contained in the disk bounded by the other, since both of the
disks bounded by the loops contain $N$. Up to isotopy in $\R^2-N$
and to a possible exchange between $C'$ and $C''$, they are as in
Figure \ref{Fig100}. The other edges of $\Gamma'$, if any, must be
contained in the annulus bounded by the two loops. So, up to an
isotopy of $\R^2-N$, which can be chosen as the identity outside
$C''$, they are as in Figure \ref{Fig4}. Of course, the same
configuration of these edges holds when $\A'=\emptyset$.

\begin{figure}[h]
\begin{center}
\includegraphics*[totalheight=4.5cm]{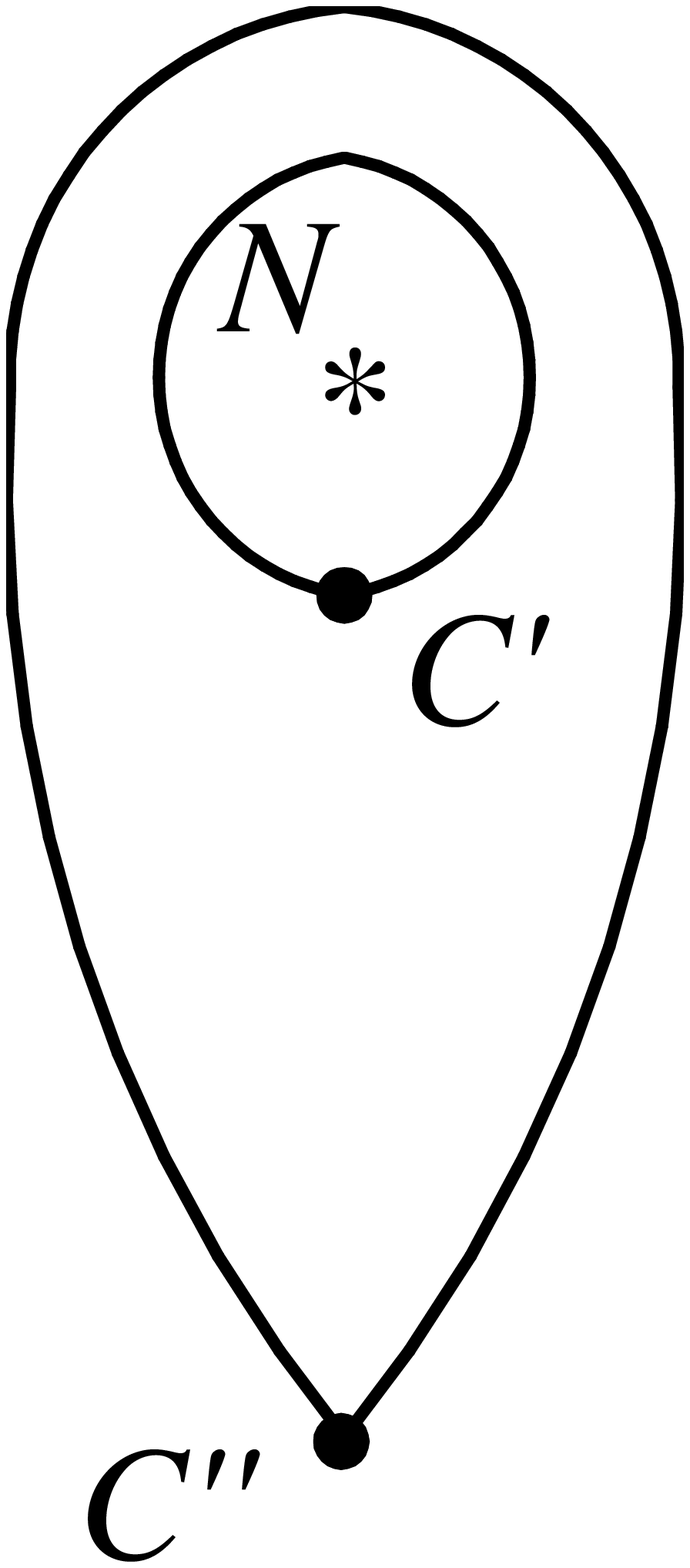}
\end{center}
\caption{} \label{Fig100}
\end{figure}

So $H_{\ps}$ is the diagram depicted in Figure \ref{Fig26}, where
$a,b,c$ are the cardinalities of the isotopy classes.
\end{proof}

\begin{figure}[h]
\begin{center}
\includegraphics*[totalheight=4.5cm]{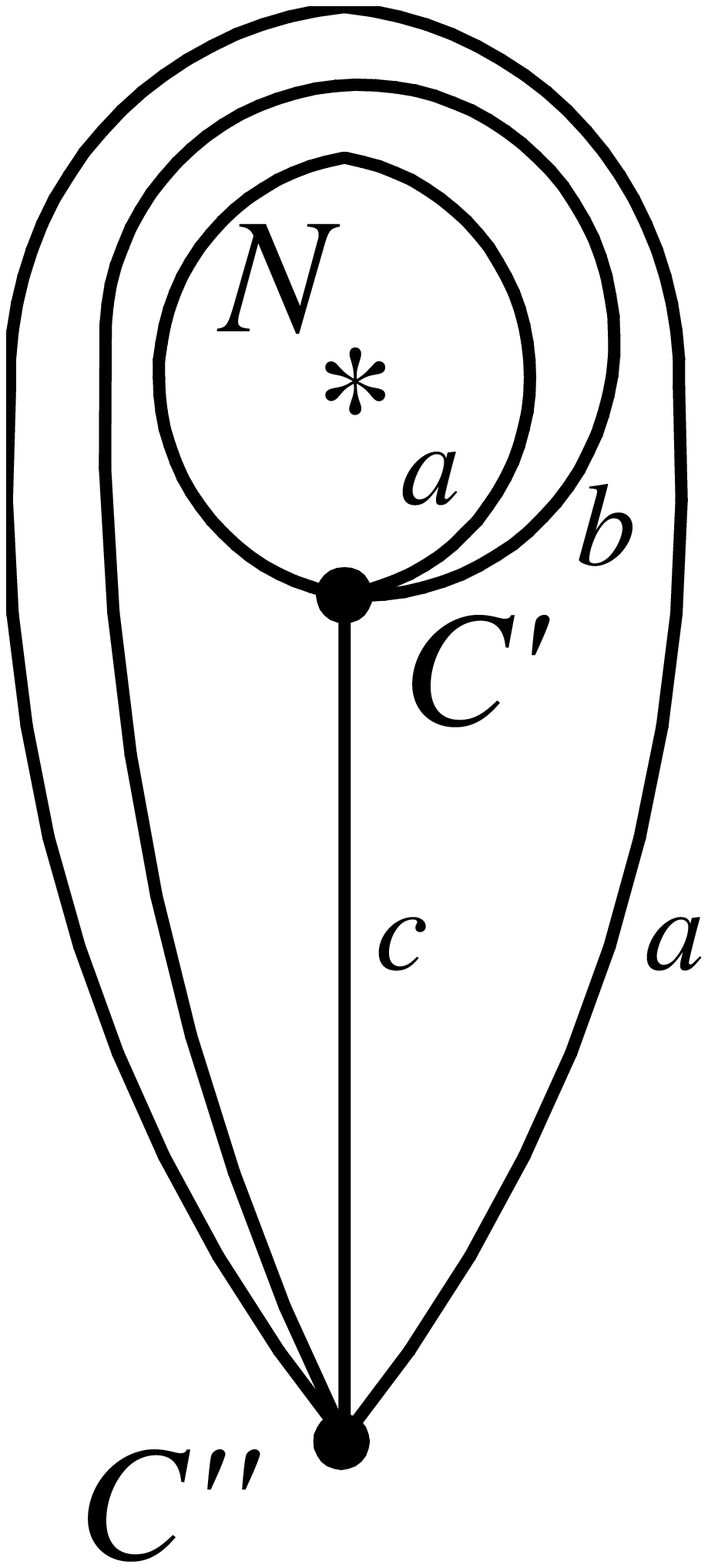}
\end{center}
\caption{} \label{Fig4}
\end{figure}

By Theorem \ref{ultimo} and Proposition \ref{strongly-cyclic} we
have:

\begin{corollary}  The class of Dunwoody manifolds coincides with the
class of strongly-cyclic branched coverings of $(1,1)$-knots.
\end{corollary}


\section{$\mathbf{(1,1)}$-knots parametrization}

As a consequence of the proof of Theorem \ref{ultimo}, any
$(1,1)$-knot $K$, with the sole exception of the core knot
$\{P\}\times\S^1\subset\S^2\times\S^1$ (which admits no
strongly-cyclic branched coverings), has a $(1,1)$-decomposition
which can be represented by an admissible Dunwoody diagram
$D(a,b,c,1,r,0)$, for suitable integers $a,b,c\ge 0$ and $r$. In
this case, we set $K=K(a,b,c,r)$, and we have that the Dunwoody
manifold $M(a,b,c,n,r,s)$ is an $n$-fold strongly-cyclic branched
covering of the lens space $M(a,b,c,1,r,0)$ (possibly homeomorphic
to $\S^3$), branched over the $(1,1)$-knot $K(a,b,c,r)$.

\medskip

\noindent {\bf Examples.} By \cite[Theorem 8]{GM}, the two-bridge
knot with Schubert parametrization $(2a+1,2r)$ is the $(1,1)$-knot
$K(a,0,1,r)$. The trivial knot in $\S^2\times\S^1$ is $K(0,0,0,0)$
and the trivial knot in $L(p,q)$ (including $L(1,0)\cong\S^3$) is
$K(0,0,p,q)$.

\medskip

Note that a different parametrization of $(1,1)$-knots, which
involves four parameters for the knot and two additional
parameters for the ambient space, can be found in \cite{CK}.

\begin{figure}[htb]
\begin{center}
\includegraphics*[totalheight=7cm]{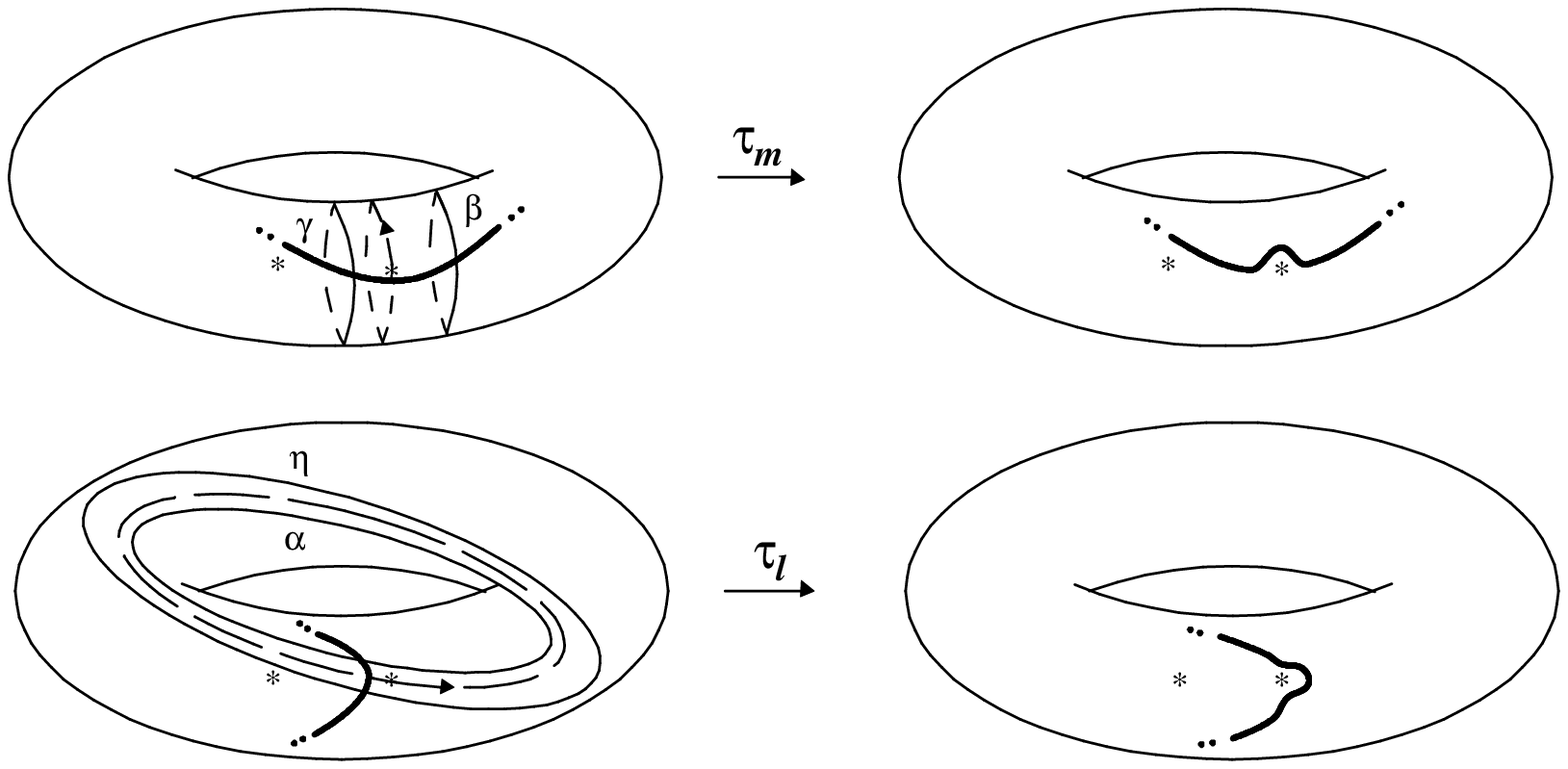}
\end{center}
\caption{Action of $\tau_m$ and $\tau_l$.} \label{Fig5}
\end{figure}

Now we describe an algorithm that gives the parametrization
$K(a,b,c,r)$  of all torus knots in $\SS^3$.

Given a closed simple curve $\delta\in\partial H$, denote by
$\t_{\delta}\in PMCG_2(\partial H)$ the right-hand Dehn twist
along $\delta$. Moreover, let $\tau_m=\t_{\b}\t_{\bb}^{-1}$ and
$\tau_l=\t_{\aa}\t_{\a}^{-1}$, where $\b,\bb,\a,\aa$ are the
curves depicted in Figure \ref{Fig5}. The effect of $\tau_m$ and
$\tau_l$ is to slide one puncture, for example $N$, along the
dashed curves depicted in Figure \ref{Fig5}, i.e. along a meridian
and a longitude of the torus, respectively.

As shown in \cite{CM2}, for every $1<k<h$, the torus knot ${\bf
t}(k,h)\subset\S^3$ is the $(1,1)$-knot $K_{\f}$ with:
\begin{equation} \label{torus}
\f=\prod_{j=0}^{h-1}(\tau_l^{-1}\tau_m^{\varepsilon_{h-j}})\t_{\b}\t_{\a}t_{\b},
\end{equation}
where\footnote{$\lfloor x \rfloor$ denotes the integral part of
$x$.} $\varepsilon_{h-j}=\lfloor (j+1)k/h \rfloor-\lfloor (j+2)k/h
\rfloor$. Since $k<h$, we have $\varepsilon_{h-j}\in\{-1,0\}$, for
all $j$.

In order to find the parameters $a,b,c,r$ for ${\bf t}(k,h)$, it
is enough to illustrate how the Heegaard diagram $D(0,0,0,1,0,0)$
is modified by the initial application of $\t_{\b}\t_{\a}t_{\b}$
and by the successive applications of the elements $\tau_l^{-1}$
and $\tau_l^{-1}\tau_m^{-1}$ composing $\f$, according to
(\ref{torus}). In this way we construct a Heegaard diagram
$D(a,b,c,1,r,0)$ representing ${\bf t}(k,h)$.

\begin{figure}
\begin{center}
\includegraphics*[totalheight=18cm]{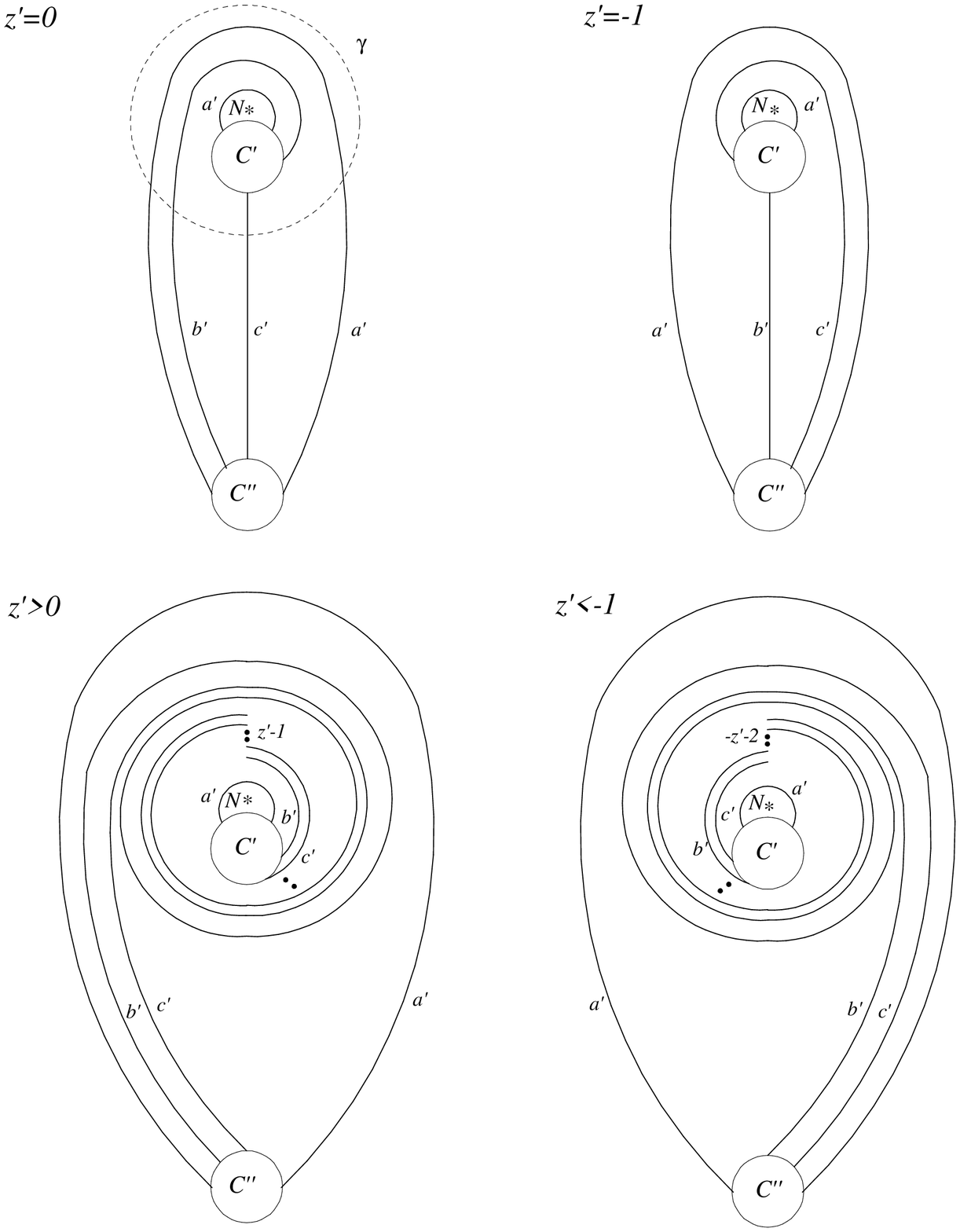}
\end{center}
\caption{The Heegaard diagram $D_{z'}(a',b',c',1,r',0)$.}
\label{Fig29}
\end{figure}

Actually, during the process, the Heegaard diagrams involved at
each step are diagrams which can be obtained by performing a
certain number $z'\in\Z$ of Dehn twists along the curve $\gamma$
to a standard Dunwoody diagram $D(a',b',c',1,r',0)$  (see Figure
\ref{Fig29}). We will call this diagram $D_{z'}(a',b',c',1,r',0)$.
These types of diagrams are depicted in Figure \ref{Fig29}, where
an arc labelled $k$ denotes $k$ parallel arcs. Obviously,
$D_0(a',b',c',1,r',0)=D(a',b',c',1,r',0)$.

Observe that, at the end of the process, we can reduce $z'$ to
zero, since $K_{t_{\gamma}\psi}$ and $K_{\psi}$ are equivalent
knots.

\begin{proposition} \label{ultima} Let $\tr(k,h)\subset\SS^3$ be a torus knot
and  $\f$ be its representation described in (\ref{torus}). Then
$\tr(k,h)=K(a,b,c,r)$ where $(a,b,c,r)=(a_h,b_h,c_h,r_h)$ is the
final step of the following algorithm, applied for
$i=h-j=1,\ldots,h$:

\begin{itemize}
\item[--] $(a_0,b_0,c_0,r_0)=(0,0,1,0)$ and $z_0=0$; \item[--] for
$i=1,\ldots,h$:
$$\left\{\begin{array}{l}a_{i}=a_{i-1}+v\\
b_{i}= r_{i-1}-2w-ud\\
c_{i}= d-b_{i}\\
r_{i}= a_{i-1}+v+w\\
z_{i}= u -\varepsilon_i
\end{array}\right.$$
where: $$ w=\begin{cases}
    a_{i-1}+b_{i-1}+c_{i-1}& \text{if $z_{i-1}<-1-\varepsilon_i$}\\
    a_{i-1}+c_{i-1}& \text{if $z_{i-1}=-1-\varepsilon_i$}\\
     a_{i-1}  & \text{if $z_{i-1}>-1-\varepsilon_i$}
      \end{cases},$$
      $$ v=\begin{cases}
    -(b_{i-1}+c_{i-1})(z_{i-1}+1+\varepsilon_i)-b_{i-1}& \text{if $z_{i-1}<-1-\varepsilon_i$} \\
    0& \text{if $z_{i-1}=-1-\varepsilon_i$}\\
    (b_{i-1}+c_{i-1})(z_{i-1}+1+\varepsilon_i)-c_{i-1}& \text{if $z_{i-1}>-1-\varepsilon_i$}
  \end{cases},$$
 and $u=\lfloor(r_{i-1}-2w)/d\rfloor$, with $d=2a_{i-1}+b_{i-1}+c_{i-1}$.
\end{itemize}
\end{proposition}

\begin{figure}
\begin{center}
\includegraphics*[totalheight=7.5cm]{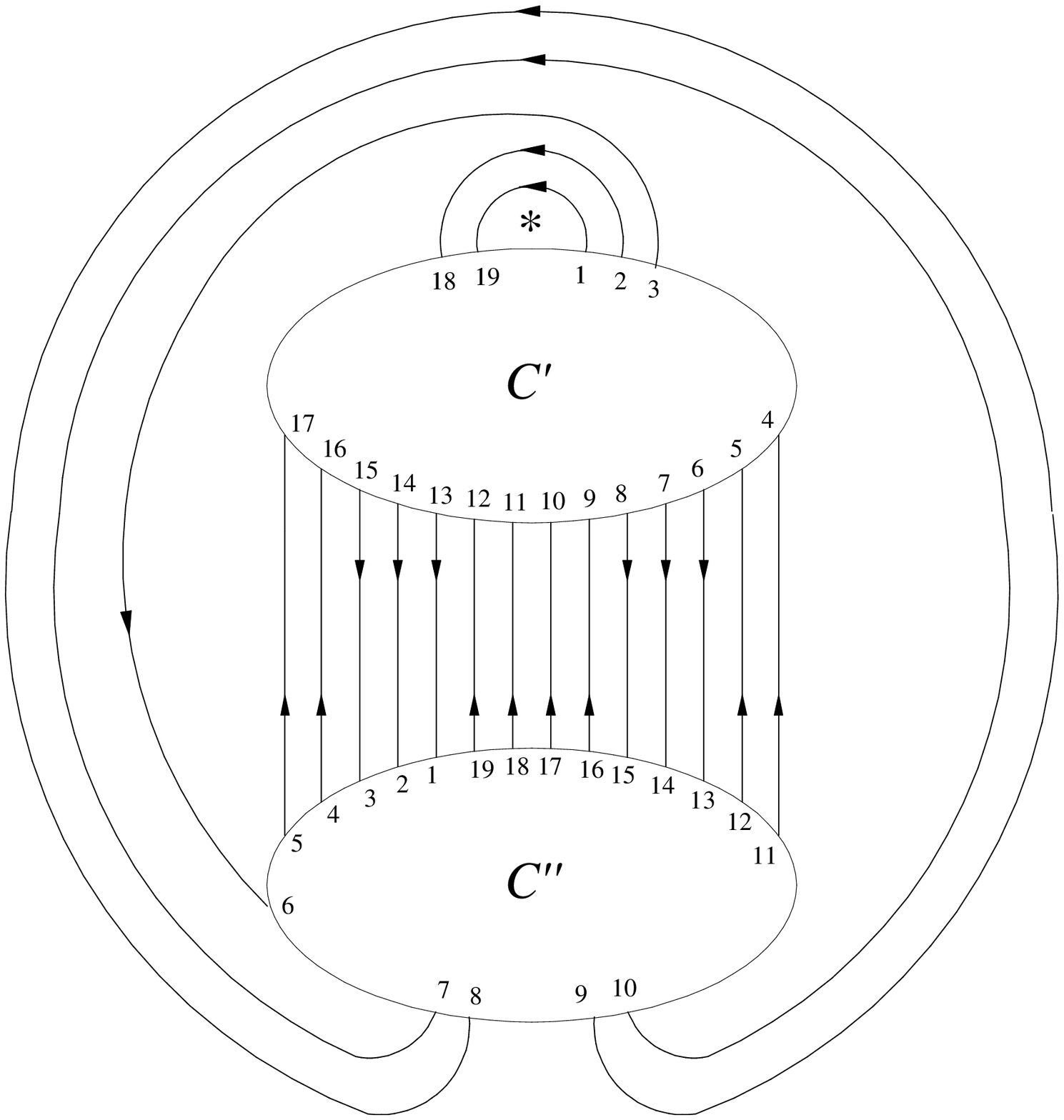}
\end{center}
\caption{$D(2,1,14,1,11,0)$.} \label{Fig39}
\end{figure}

The proof of Proposition \ref{ultima} will be given at the end of
this section. Now we give some examples and applications.

\begin{remark} \label{rives} Given an admissible Dunwoody diagram
$D(a,b,c,1,r,0)$, with $a+b+c>0$, we  fix an orientation on the
arcs of $\E$  that induces  an orientation on the corresponding
curve of the Heegaard diagram in such a way that the vertex on
$C'$ labelled 1 is the first endpoint of the corresponding edge.
Let $p_{a,b,c,r}$ be the number of arcs of $\B$ oriented from $C'$
to $C''$ minus the number of arcs oriented from $C''$ to $C'$, and
let $q_{a,b,c,r}$ be the number of arcs of $\E$ oriented from
right to left minus the number of arcs oriented from left to right
(see \cite[p. 385]{GM}). If $K(a,b,c,r)$ is a $(1,1)$-knot in
$\SS^3$, then the $n$-fold cyclic branched covering of
$K(a,b,c,r)$ is the Dunwoody manifold $M(a,b,c,r,n,s)$, where
$s=-p_{a,b,c,r}q_{a,b,c,r}$. In fact,  by Proposition
\ref{strongly-cyclic}, there exists a unique $s$ (mod $n$) such
that $M(a,b,c,n,r,s)$ is the $n$-fold cyclic covering of
$M(a,b,c,1,r,0)\cong \SS^3$, branched over $K(a,b,c,r)$. Moreover,
by \cite{GM}, $s$ must satisfy the condition
$q_{a,b,c,r}+sp_{a,b,c,r}\equiv 0$ (mod $n$) and we have
$p_{a,b,c,r}=\pm 1$.
\end{remark}


\noindent {\bf Example.} Let us consider $\tr(5,8)$. By
(\ref{torus}), a representation of $\tr(5,8)$ is given by
\hbox{$\f=\tau_l^{-1}\tau_m^{-1}\tau_l^{-1}(\tau_l^{-1}(\tau_m^{-1}\tau_l^{-1})^2)^2\t_{\b}\t_{\a}t_{\b}$.}
Then, by Proposition \ref{ultima}, we have
$\tr(5,8)=K(2,1,14,11)$.  Moreover, from the diagram
$D(2,1,14,1,11,0)$ depicted in Figure \ref{Fig39}, we get
$p_{2,1,14,11}=-1$ and $q_{2,1,14,11}=5$. So, by Remark
\ref{rives}, the $n$-fold cyclic branched covering of $\tr(5,8)$
is the Dunwoody manifold $M(2,1,14,n,11,5)$, for all $n>1$.

\medskip

As an application, we explicitly determine the parametrization of
\hbox{$\tr(k,ck+1)$} as well as the Dunwoody representation of its
cyclic branched coverings.
\begin{corollary} \label{fine}
For every $c>0$ and $k>1$,  the torus knot $\tr(k,ck+1)$ is
$K(1,k-2,2kc-2c-k+1,k)$. Moreover, the  $n$-fold cyclic branched
covering of $\tr(k,ck+1)$ is the Dunwoody manifold
$M(1,k-2,2kc-2c-k+1,n,k,k)$, for all $n>1$.
\end{corollary}
\begin{proof}  By (\ref{torus}), $\tr(k,ck+1)$
is  represented by
$\f=(\tau_l^{-c}\tau_m^{-1})^k\tau_l^{-1}\t_{\b}\t_{\a}\t_{\b}$.
Applying  Proposition \ref{ultima} and   Remark \ref{rives} we get
the statement.
\end{proof}

Observe that Corollary \ref{fine} agrees with the result obtained
in \cite{AGM} with different techniques.

\medskip

\begin{figure}
\begin{center}
\includegraphics*[totalheight=6cm]{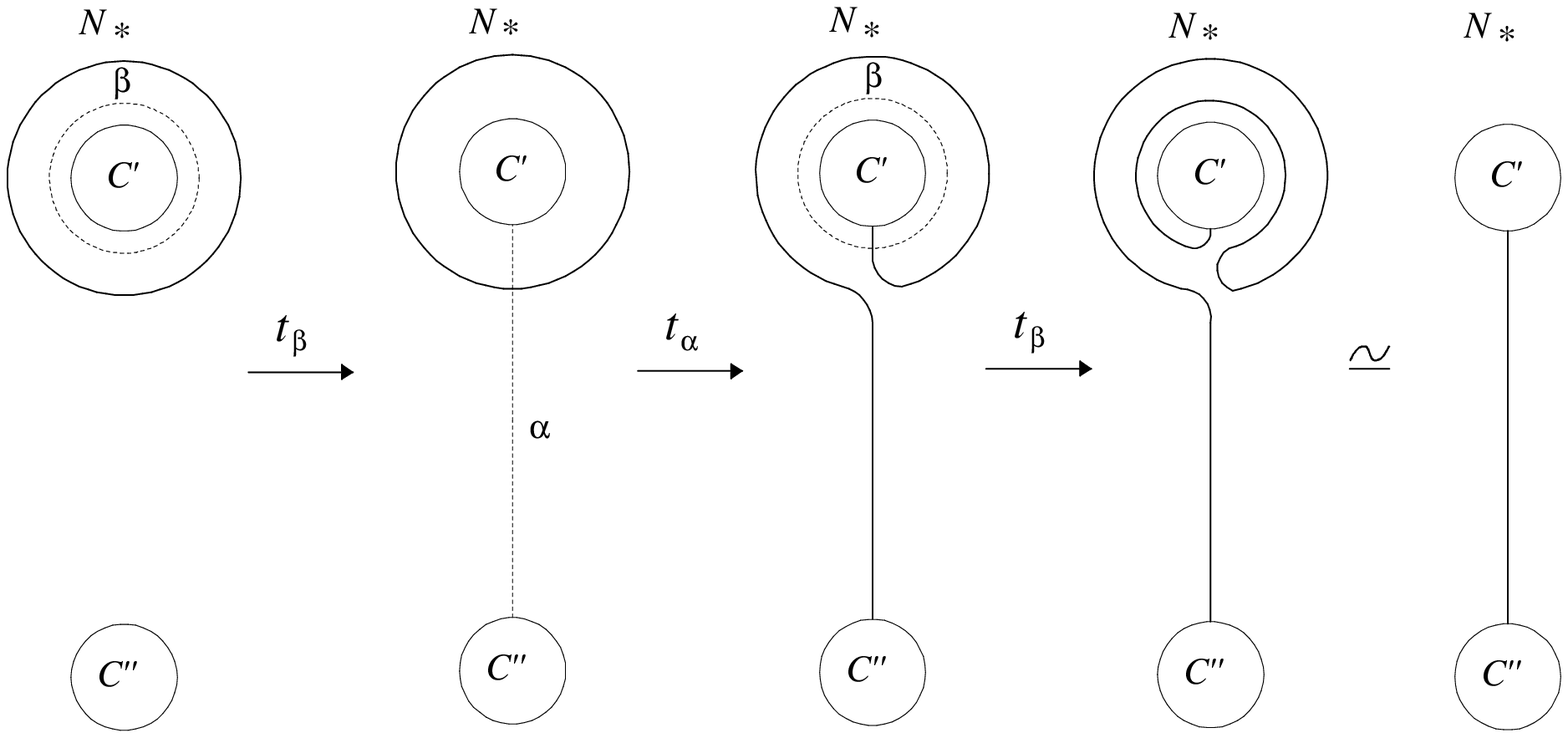}
\end{center}
\caption{Action of $\t_{\b}\t_{\a}\t_{\b}$ on $D(0,0,0,1,0,0)$.}
\label{Fig33}
\end{figure}

\noindent \textit{Proof of  Proposition \ref{ultima}.} As shown in
Figure \ref{Fig33}, the application of $\t_{\b}\t_{\a}\t_{\b}$ to
$D(0,0,0,1,0,0)$ gives the diagram $D(0,0,1,1,0,0)$.

In order to simplify the notations in the figures, we set
$(a_{i-1},b_{i-1},c_{i-1},r_{i-1})=(a',b',c',r')$ and
$z_{i-1}=z'$. To obtain the parameters $a,b,c$ and $r$, we
consider the application of $\tau_l^{-1}\tau_m^{\varepsilon_i}$ to
$D_{z'}(a',b',c',1,r',0)$.

\begin{figure}
\begin{center}
\includegraphics*[totalheight=18cm]{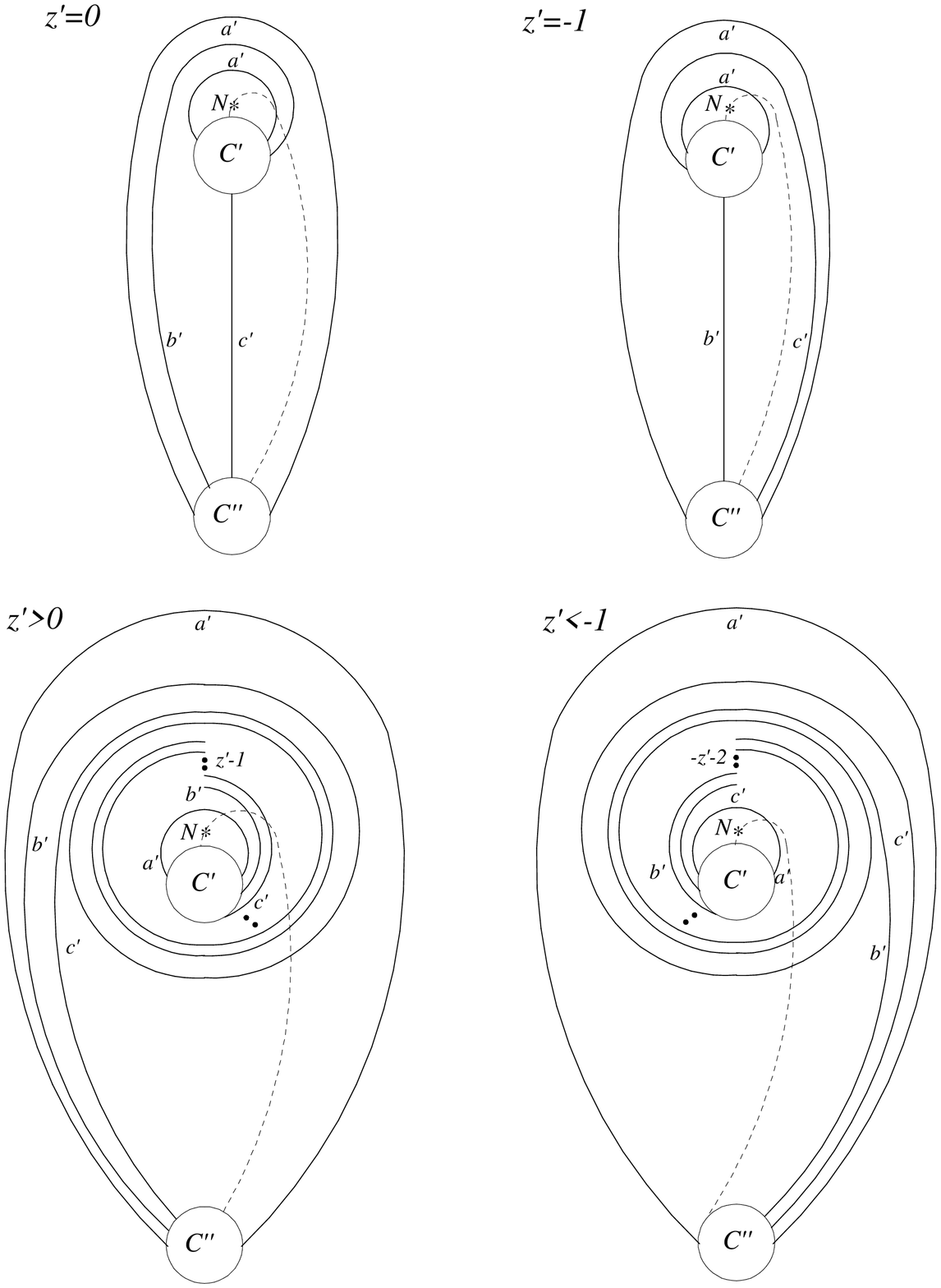}
\end{center}
\caption{The parameter  $w$.} \label{Fig35}
\end{figure}

Let us  first consider the case  $\varepsilon_i=0$. We recall that
the effect of $\tau_l^{-1}$ is to slide $N$ along the longitude of
the torus, illustrated  by the dashed line in Figure \ref{Fig5},
in the opposite direction to the arrow. This curve will always be
represented on a Heegaard diagram by a dashed arc connecting an
internal point of the arc on $C'$, with endpoints labelled $d$ and
$1$ (according to the orientation), with the corresponding point
on $C''$. The number of intersections of the longitude with the
arcs of a given diagram depends on $r'$. Let $w$ be the value of
$r'$ such that the number of these intersections is minimal. Then,
as illustrated in Figure \ref{Fig35}, we have:
$$w=\begin{cases}
    a'+b'+c'& \text{if $z'<-1$} \\
    a'+c' & \text{if $z'=-1$}\\
     a'  & \text{if $z'>-1$}
\end{cases}.$$

In this figure, and in the following ones, an arc labelled $f$
denotes $f$ parallel arcs, and we take the convention that a label
of a vertex is the label corresponding to the endpoint of the
first of the $f$ parallel arcs.

\begin{figure}
\begin{center}
\includegraphics*[totalheight=17cm]{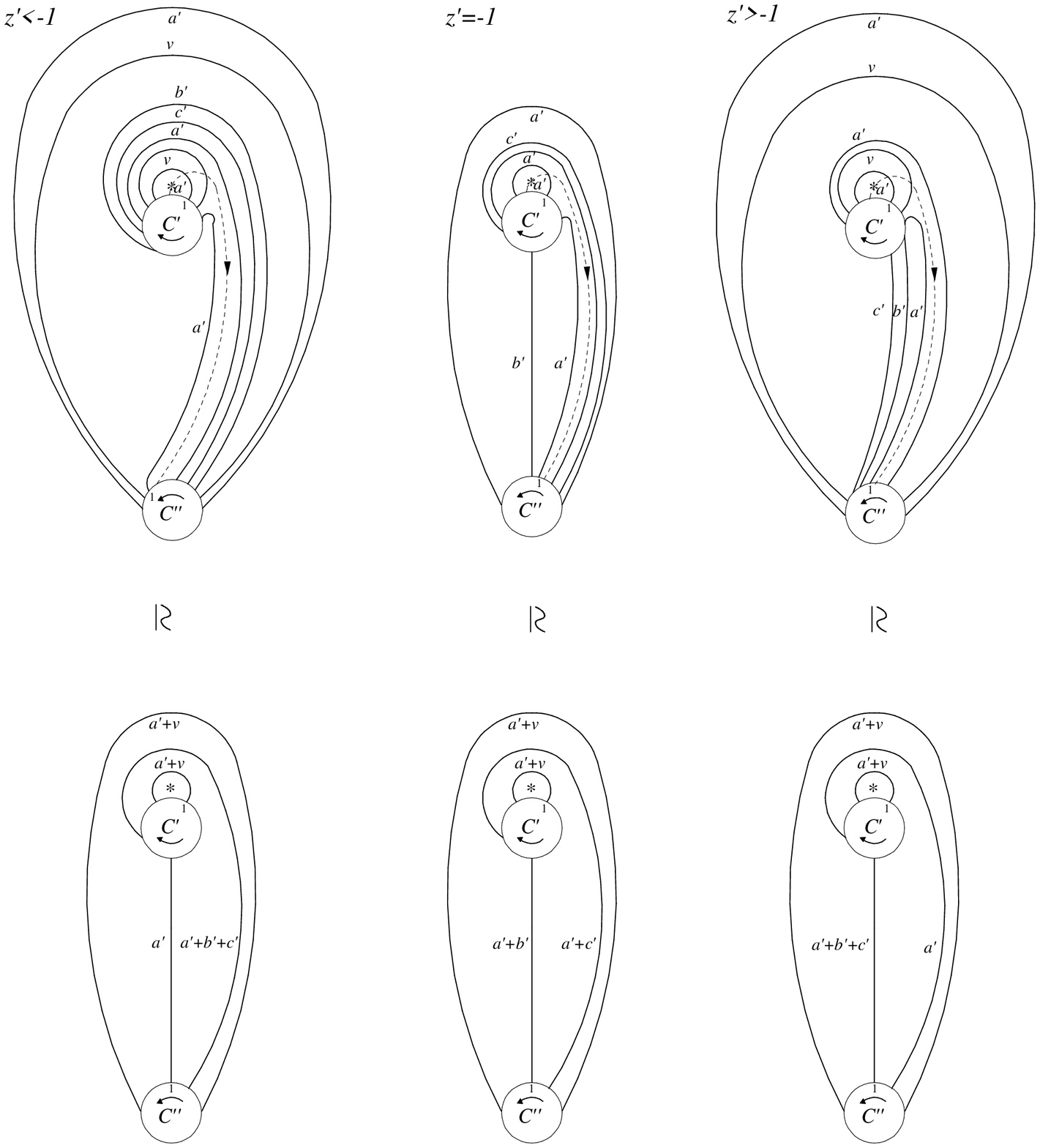}
\end{center}
\caption{Action of  $\tau_l^{-1}$ on $D_{z'}(a',b',c',1,r',0)$ for
$r'=w$.} \label{Fig36}
\end{figure}

First of all, we consider the case $r'=w$. In this case the
longitude has  $a'+v$ intersections, and the action of
$\tau_l^{-1}$ is illustrated in Figure \ref{Fig36}.
We obtain $(a_i,b_i,c_i,r_i)=(a'+v,d-w,w,a'+v+w)$ and $z_i=-1$,
which is the same result of the statement when $r'=w\neq 0$ (in
this case $-d\leq r'-2w=-w<0$ and so $u=-1$).  If $r'=w=0$, we
have $u=0$, and therefore the statement gives
$(a_i,b_i,c_i,r_i)=(a'+v,0,d,a'+v)$ and $z_i=0$; but it is easy to
check that $D_0(a'+v,0,d,1,a'+v,0)=D_{-1}(a'+v,d,0,1,a'+v,0)$.

\begin{figure}
\begin{center}
\includegraphics*[totalheight=8cm]{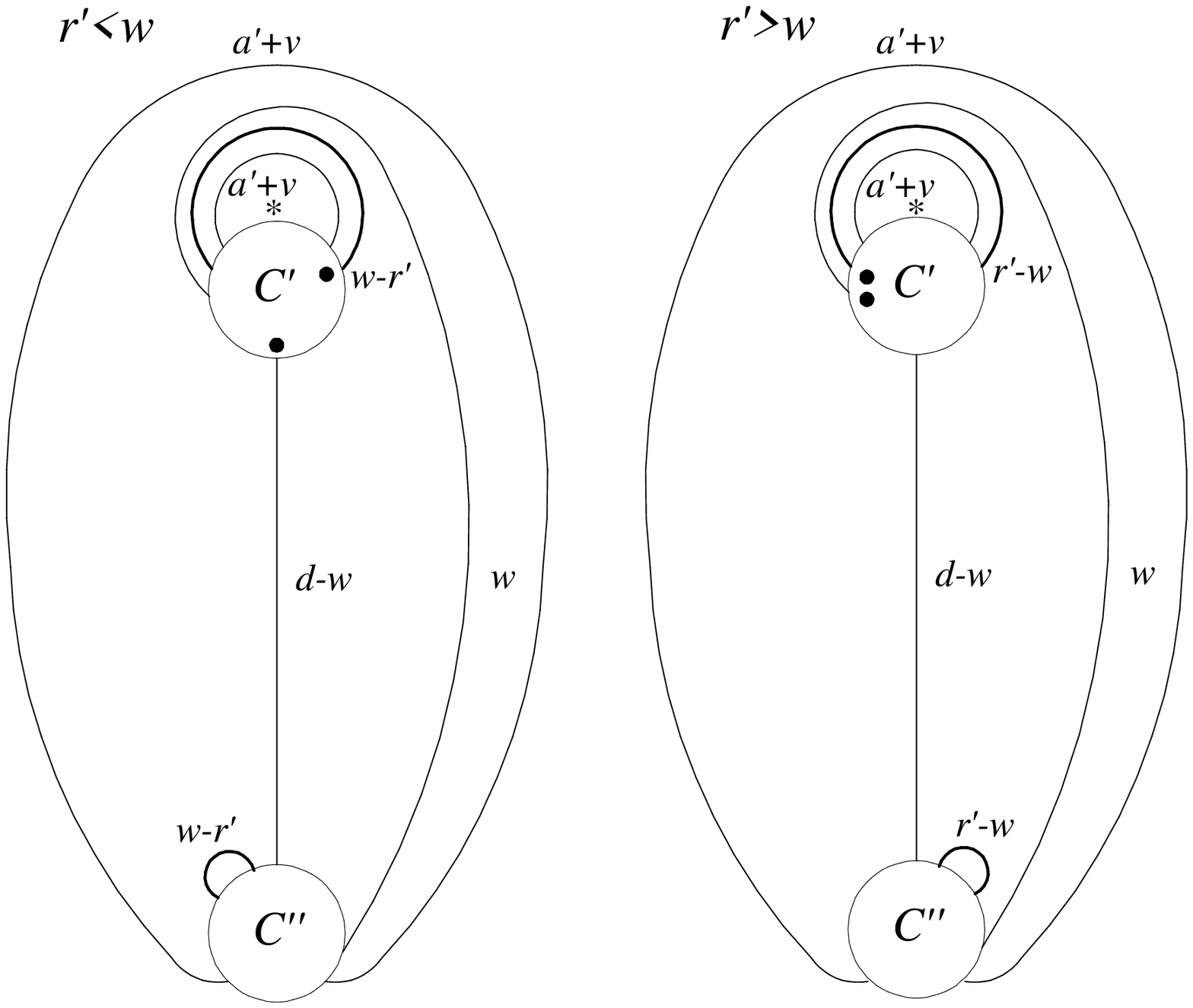}
\end{center}
\caption{Action of $\tau_l^{-1}$ on $D_{z'}(a',b',c',1,r',0)$ for
$r'<w$ and $r'>w$.} \label{Fig37}
\end{figure}

\begin{figure}
\begin{center}
\includegraphics*[totalheight=8cm]{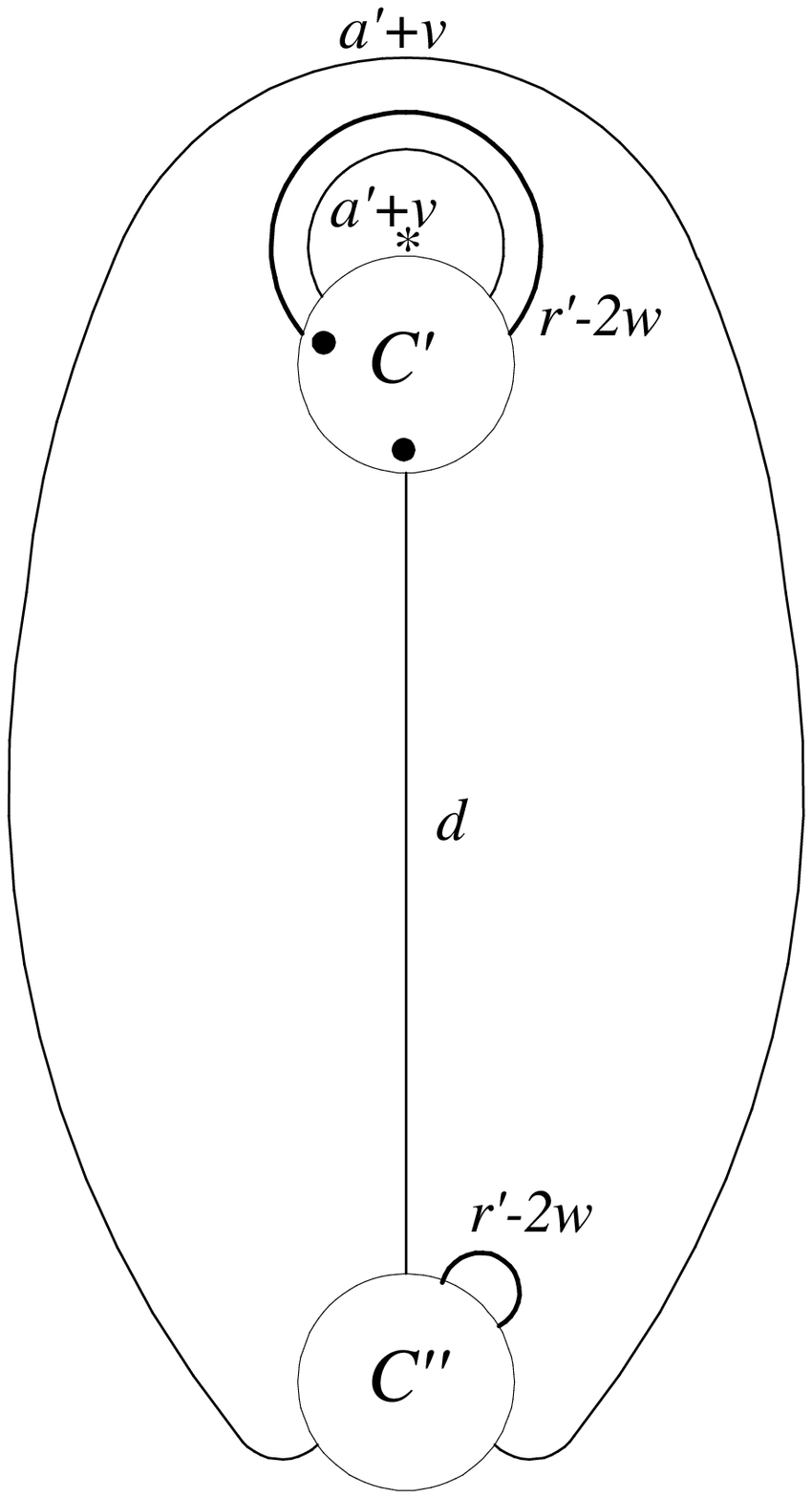}
\end{center}
\caption{Action of $\tau_l^{-1}$ in the case $r'-w\geq w$.}
\label{Fig38}
\end{figure}

When $r'>w$ or $r'<w$, the result of the application of
$\tau_l^{-1}$ is depicted  in Figure \ref{Fig37}. In both cases,
the further $|r'-w|$ intersections determine $|r'-w|$ trivial arcs
on $C''$. The $j$-th of these arcs has endpoints on $C''$ labelled
$a'+v+d+j$ and $a'+v+d+2(r'-w)-j+1$ if $w<r'$, and labelled
$a'+v+j$ and $a'+v+2(w-r')-j+1$ if $w>r'$.  Each time we eliminate
a trivial arc $e$, we glue together the two arcs whose endpoints
on $C'$ have the same label as the endpoints of $e$ on $C''$. In
Figure \ref{Fig37}, the black points indicate which arcs are glued
together. After the elimination of all the trivial arcs, we
obtain, as above, $a_i=a'+v$ and $r_i=a'+v+w$, while the value of
the other three parameters depends on the quotient of the division
of $|r'-2w|$ by $d$. Suppose that $r'>w$, then we have two cases:
\begin{itemize}
\item[\textup{(1)}] if $r'-w<w$, we obtain $b_i=d-w+r'-w=d+r'-2w$,
\hbox{$c_i=w-(r'-w)=2w-r'$} and $z_i=-1$; \item[\textup{(2)}] if
$r'-w\geq w$, after the elimination of the first $w$ trivial arcs,
we obtain the diagram depicted  in Figure \ref{Fig38}. During the
elimination of the remaining $r'-2w$ arcs, each time we eliminate
$d$ arcs the parameter $z'$ increases by one. Therefore, if $u$ is
the integer defined by $u=\lfloor(r'-2w)/d\rfloor$, we have
$b_i=r'-2w-ud$, $c_i=(u+1)d-(r'-2w)$ and $z_i=u$.
\end{itemize}

\begin{figure}
\begin{center}
\includegraphics*[totalheight=8cm]{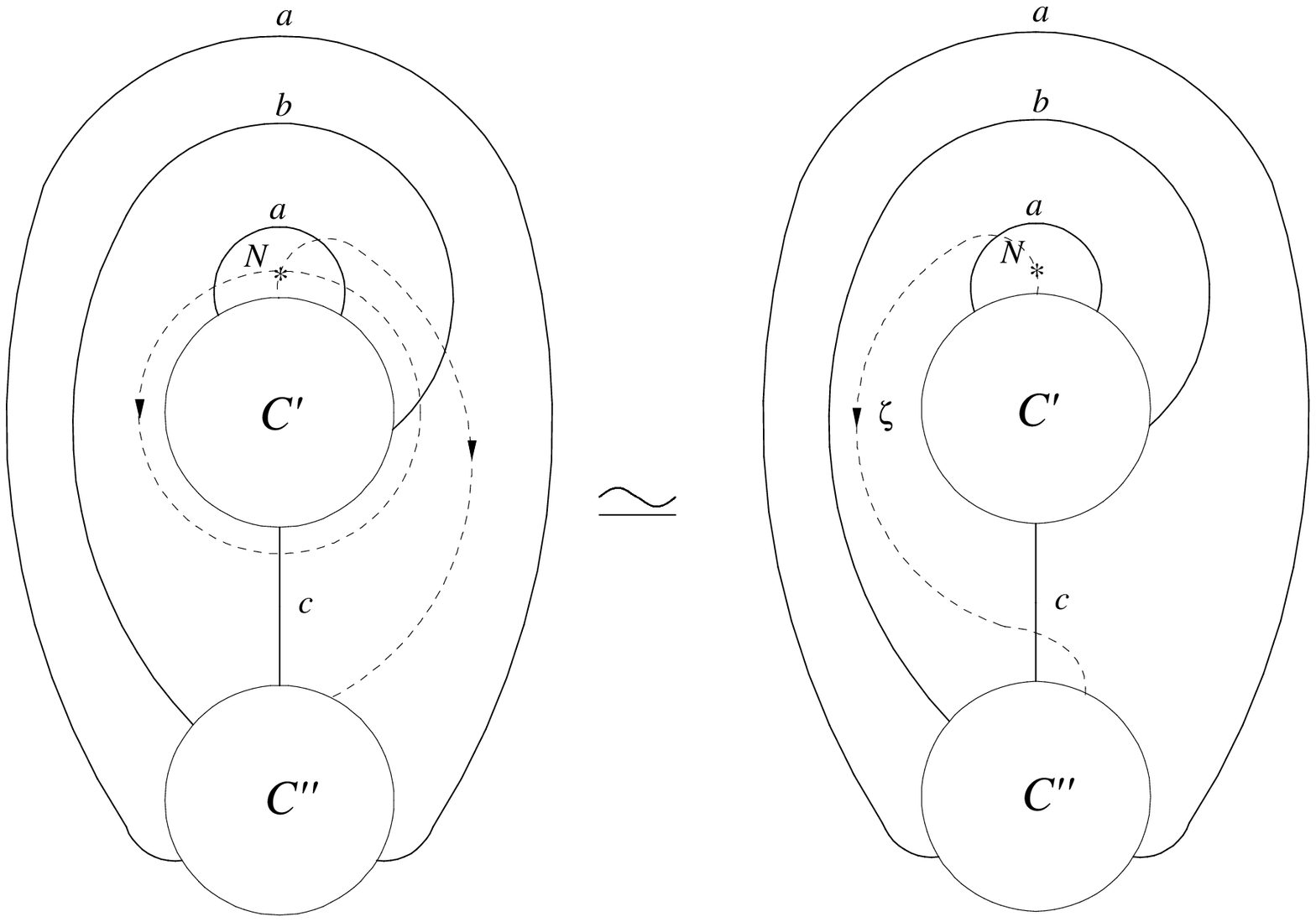}
\end{center}
\caption{Action of $\tau_l^{-1}\tau_m^{-1}$.} \label{Fig0}
\end{figure}

Analysing the case $r'<w$ in an analogous way, we complete the
case $\varepsilon_i=0$.

In the  case $\varepsilon_i=-1$ we examine the action of
$\tau_l^{-1}\tau_m^{-1}$. This can be done in a similar way as
before, since, as depicted in Figure \ref{Fig0}, the action of
$\tau_l^{-1}\tau_m^{-1}$ is equivalent to an action that moves $N$
along the longitude $\zeta$. \hspace{\stretch{1}}$\square$

\bigskip

\noindent{\bf Acknowledgements}

\noindent Work performed under the auspices of the G.N.S.A.G.A. of
I.N.d.A.M. (Italy) and the University of Bologna, funds for
selected research topics.



\vspace{15 pt} {ALESSIA CATTABRIGA, Department of Mathematics,
University of Bologna, Italy. E-mail: cattabri@dm.unibo.it}

\vspace{15 pt} {MICHELE MULAZZANI, Department of Mathematics and
C.I.R.A.M., University of Bologna, Italy. E-mail:
mulazza@dm.unibo.it}

\end{document}